\documentclass[11pt]{amsart}
\usepackage{epsfig}
\usepackage[all]{xy}
\begin{document}
\title{Sur des vari\'et\'es
toriques non projectives}
\author{Laurent BONAVERO}
\date{Mai 1999}
\maketitle
\noindent
\def\restriction{\string |}
\newcommand{\pp}{\rm ppcm}
\newcommand{\pg}{\rm pgcd}
\newcommand{\Ker}{\rm Ker}
\newcommand{\C}{{\mathbb C}}
\newcommand{\Q}{{\mathbb Q}}
\newcommand{\GL}{\rm GL}
\newcommand{\SL}{\rm SL}
\newcommand{\diag}{\rm diag}

\def\refname{R\'ef\'erences}
\def\finpreuve
{\hskip 3pt \vrule height6pt width6pt depth 0pt}

\newtheorem{theo*}{Th\'eor\`eme}
\newtheorem{prop}{Proposition}
\newtheorem{lemm}{Lemme}
\newtheorem{lemmf}{Lemme fondamental}
\newtheorem{defi}{D\'efinition}
\newtheorem{exo}{Exercice}
\newtheorem{rem}{Remarque}
\newtheorem{cor}{Corollaire}
\newcommand{\CC}{{\mathbb C}}
\newcommand{\ZZ}{{\mathbb Z}}
\newcommand{\RR}{{\mathbb R}}
\newcommand{\QQ}{{\mathbb Q}}
\newcommand{\FF}{{\mathbb F}}
\newcommand{\PP}{{\mathbb P}}
\newcommand{\codim}{\operatorname{codim}}
\newcommand{\Ho}{\operatorname{Hom}}
\newcommand{\Pic}{\operatorname{Pic}}
\newcommand{\NE}{\operatorname{NE}}
\newcommand{\Nun}{\operatorname{N}}
\newcommand{\card}{\operatorname{card}}
\newcommand{\Hilb}{\operatorname{Hilb}}
\newcommand{\mult}{\operatorname{mult}}
\newcommand{\vol}{\operatorname{vol}}
\newcommand{\divi}{\operatorname{div}}
\newcommand{\pr}{\operatorname{pr}}
\newcommand{\con}{\operatorname{cont}}
\newcommand{\ima}{\operatorname{Im}}

\newcounter{subsub}[subsection]
\def\thesubsub{\thesubsection .\arabic{subsub}}
\def\subsub#1{\addtocounter{subsub}{1}\par\vspace{3mm}
\noindent{\bf \thesubsub ~ #1 }\par\vspace{2mm}}
\def\coker{\mathop{\rm coker}\nolimits}
\def\pr{\mathop{\rm pr}\nolimits}
\def\im{\mathop{\rm Im}\nolimits}
\def\hfl#1#2{\smash{\mathop{\hbox to 12mm{\rightarrowfill}}
\limits^{\scriptstyle#1}_{\scriptstyle#2}}}
\def\vfl#1#2{\llap{$\scriptstyle #1$}\big\downarrow
\big\uparrow
\rlap{$\scriptstyle #2$}}
\def\diagram#1{\def\normalbaselines{\baselineskip=0pt
\lineskip=10pt\lineskiplimit=1pt}   \matrix{#1}}
\def\limind{\mathop{\oalign{lim\cr
\hidewidth$\longrightarrow$\hidewidth\cr}}}

\long\def\InsertFig#1 #2 #3 #4\EndFig{
\hbox{\hskip #1 mm$\vbox to #2 mm{\vfil\includegraphics{#3}}#4$}}
\long\def\LabelTeX#1 #2 #3\ELTX{\rlap{\kern#1mm\raise#2mm\hbox{#3}}}

{\let\thefootnote\relax
\footnote{%\hskip3em 
\textbf{Mots cl\'es :} vari\'et\'e torique, vari\'et\'e 
non projective, th\'eorie de Mori.
\textbf{Classification~A.M.S. :} 14M25, 14E05, 14E30. 
}}

{\bf R\'esum\'e : } Utilisant la th\'eorie de Mori des vari\'et\'es
toriques projectives due \`a M. Reid, nous \'etudions
les vari\'et\'es toriques non projectives qui 
le deviennent
apr\`es \'eclatement le long d'une courbe invariante. 

\medskip

{\bf Abstract : } With the help of Mori theory
for projective toric manifolds due to M. Reid, 
we study non projective toric manifolds
which become projective after a single blow up along an
invariant curve.

\section*{Introduction}

Dans tout ce travail, nous appelons vari\'et\'e torique
de dimension $n$ une vari\'et\'e torique compl\`ete et non
singuli\`ere,
autrement dit, nous nous int\'eressons aux compactifications lisses
\'equivariantes du tore complexe $(\CC ^*)^n$.
Une telle vari\'et\'e n'est en g\'en\'eral pas projective, mais
les exemples connus de vari\'et\'es toriques non projectives
peuvent \^etre consid\'er\'es comme isol\'es au sens 
o\`u la cause g\'eom\'etrique de non projectivit\'e n'est pas 
toujours comprise. Nous nous proposons ici d'\'etudier 
g\'eom\'etriquement les vari\'et\'es toniques 
non projectives les plus simples. 
Pour cela, rappelons qu'il est bien connu (voir par exemple \cite{Mor96}) 
que toute vari\'et\'e torique
devient projective apr\`es une suite finie d'\'ecla\-te\-ments
le long de sous-vari\'et\'es lisses toriques~; autrement
dit, la version torique du th\'eor\`eme de Moishezon \cite{Moi67}
est v\'erifi\'ee. 
{\em Le but de ce travail est de comprendre les
vari\'et\'es toriques non projectives les plus simples 
\`a la lumi\`ere du th\'eor\`eme de Moishezon~: celles qui deviennent
projectives apr\`es \'eclatement le long d'une courbe torique }
($X$ \'etant lisse, toute courbe torique de $X$ est 
rationnelle non singuli\`ere).
Rappelons aussi que si $X$ est une vari\'et\'e complexe compacte
et si $x$ est un point de $X$, alors
$X$ est projective si et seulement si $B_x(X)$
est projective. 

\medskip

\noindent {\bf Notation. } {\em Si $X$ est 
une vari\'et\'e lisse et $Z$ une sous-vari\'et\'e lisse
de $X$, on note $B_Z(X)$ la vari\'et\'e
obtenue en \'eclatant $X$ le long de $Z$. 
Rappelons que si $E$ est le diviseur exceptionnel de
l'\'eclatement, $E$ est isomorphe \`a $\PP (N_{Z/X})$
o\`u $N_{Z/X}$ est le fibr\'e normal de $Z$ dans $X$.
Si $X$ est une vari\'et\'e torique et $Z$ est une sous-vari\'et\'e
torique, alors $B_Z(X)$ est torique 
et l'\'eclatement~$\pi : B_Z(X) \to X$
est \'equivariant.}

\medskip

Notre travail s'articule de la fa\c con suivante~:
nous rappelons dans un premier temps les notations 
classiques en g\'eom\'etrie torique et l'\'enonc\'e principal
de la th\'eorie de Mori des vari\'et\'es toriques projectives
d'apr\`es Miles Reid \cite{Rei83}. 
On en d\'eduit un crit\`ere de projectivit\'e
pour les vari\'et\'es toriques qui sont projectives
apr\`es un \'eclatement. 
Dans un deuxi\`eme temps, nous montrons que si $X$ est une
vari\'et\'e torique non projective
contenant une courbe $C$ torique de sorte que 
$\tilde X~:=B_C(X)$ soit projective, alors $C$ est rigide
dans $X$.
Dans les troisi\`eme et quatri\`eme parties, si $X$ est une
vari\'et\'e torique non projective
contenant une courbe $C$ torique de sorte que 
$\tilde X~:=B_C(X)$ soit projective, 
nous classifions les contractions 
{\em Mori extr\'emales} 
sur $\tilde X$ dont le lieu exceptionnel rencontre
le diviseur exceptionnel de l'\'eclatement $\tilde X \to X$. Nous montrons 
en particulier que toutes sont des 
\'eclatements le long de centres lisses dans une
vari\'et\'e lisse et que seuls
trois ph\'enom\`enes peuvent se produire dans le
cas o\`u $\tilde X$ poss\`ede une telle contraction~: soit 
$(X,C)$ poss\`ede une {\em r\'eduction triviale}, 
soit $X$ est obtenue {\em via} un {\em flip 
interdit } d'une vari\'et\'e torique projective, soit $X$
est obtenue {\em via} une {\em transformation
\'el\'ementaire} (dite ``de Maruyama'')
d'une vari\'et\'e torique projective, 
ces termes \'etant introduits et illustr\'es au {\S}3.
Nous donnons aussi deux cons\'equences pour les vari\'et\'es toriques 
non projectives
devenant de Fano apr\`es \'eclatement le long d'une courbe, en montrant
en particulier que ceci ne se produit pas en dimension trois.
Ce dernier r\'esultat n'utilise pas
la classification des vari\'et\'es toriques de Fano
de dimension trois. 
Dans la derni\`ere partie, nous reprenons et d\'etaillons une construction
d'Ewald, permettant de donner des exemples de vari\'et\'es
toriques non projectives pour toutes les situations 
\'etudi\'ees pr\'ec\'edemment, ces vari\'et\'es poss\'edant 
de plus un petit groupe de Picard. 

\medskip

\noindent {\em Remerciements~:} \`a Michel Brion et Laurent
Manivel pour de nombreuses discussions et
leurs lectures critiques de versions pr\'eliminaires, 
\`a Alexis Marin qui m'a permis
d'interpr\'eter g\'eom\'etriquement la construction d'Ewald
et \`a Chris Peters qui m'a indiqu\'e les travaux de Maruyama.

\section{Pr\'eliminaires de g\'eom\'etrie torique}

\subsection{Premi\`eres d\'efinitions et notations}

Les r\'ef\'erences usuelles de g\'eom\'etrie torique sont
\cite{Ewa96}, \cite{Ful93} et \cite{Oda88}.

Si $M$ est un r\'eseau de rang $n$, une vari\'et\'e torique $X_{\Delta}$
de dimension $n$ est d\'efinie par la donn\'ee d'un \'eventail
$\Delta$, subdivision de l'espace vectoriel dual
$N_{\QQ}:= \Ho (M,\ZZ)\otimes_{\ZZ} \QQ$ par des c\^ones rationnels
poly\'edraux.
Rappelons qu'il y a une correspondance bijective entre les 
orbites (du tore $\Ho (M, \CC ^*)$) de codimension $r$ dans $X_{\Delta}$ 
et les c\^ones de dimension $r$ de $\Delta$. Comme annonc\'e
dans l'introduction, toutes les vari\'et\'es toriques $X_{\Delta}$
seront suppos\'ees lisses (ce qui signifie que tous 
les c\^ones maximaux de $\Delta$ sont simpliciaux, engendr\'es
par une base de $N:=\Ho (M,\ZZ)$) et compl\`etes (ce qui signifie
que le support de $\Delta$ est $ N_{\QQ}$). 

Les g\'en\'erateurs dans $N$ des faces de dimension un de $\Delta$ seront 
not\'es $e_1,\ldots,e_{\rho+n}$~; ils cor\-res\-pondent
aux diviseurs irr\'eductibles invariants de $X_{\Delta}$ 
et $\rho$ est le nombre de Picard de~$X_{\Delta}$.   

\subsection{\'Eclatements le long de courbes toriques}
  
Soit $X$ une vari\'et\'e torique d'\'eventail $\Delta$
et $C$ une courbe torique de $X$. Quitte \`a renum\'eroter 
les $e_i$, $C$ est donn\'ee par le c\^one de
dimension $n-1$ (appel\'e aussi ``mur'') $<e_1,\ldots,e_{n-1}>$.
Ce mur est face de deux c\^ones maximaux $<e_1,\ldots,e_{n-1},e_n>$
et $<e_1,\ldots,e_{n-1},e_{n+1}>$
correspondant
aux deux points fixes du tore sur $C$.

La vari\'et\'e $\tilde X = B_C(X)$ est d\'efinie par l'\'eventail
$\tilde {\Delta}$ construit de la fa\c con
suivante~: on note $e=e_1+\cdots+e_{n-1}$,
on enl\`eve \`a $\Delta$ les deux c\^ones maximaux 
$<e_1,\ldots,e_{n-1},e_n>$
et $<e_1,\ldots,e_{n-1},e_{n+1}>$ (ainsi que leurs faces) et on les
remplace par les $2n-2$ c\^ones maximaux de la forme
$$<e_1,\ldots,\hat{e_j},\ldots,e_{n-1},e,e_n> \, \mbox{ et } \,
<e_1,\ldots,\hat{e_j},\ldots,e_{n-1},e,e_{n+1}>, \,\,  1\leq j \leq n-1$$
(avec toutes leurs faces) (la notation $\hat{e_j}$ signifie qu'on 
enl\`eve $e_j$) - voir la figure~1.

Les courbes toriques irr\'eductibles sur le diviseur exceptionnel
$E$ de l'\'eclatement $\pi~:\tilde X \to X$
sont les $n-1$ courbes
$\tilde{C}_j := <e_1,\ldots,\hat{e_j},\ldots,e_{n-1},e>$, $1\leq j \leq n-1$
qui se projettent isomorphiquement par $\pi$ sur $C$
et pour $1\leq i < j \leq n-1$ les courbes 
$$<e_1,\ldots,\hat{e_i},\ldots,\hat{e_j},\ldots,e_{n-1},e,e_n> \, \mbox{et}
\,  
<e_1,\ldots,\hat{e_i},\ldots,\hat{e_j},\ldots,e_{n-1},e,e_{n+1}>$$
contract\'ees par $\pi$.

\begin{figure}[hbtp]
  \begin{center}
    \leavevmode
    \input{fig1.pstex_t}
    \caption{}
    %\label{fig:}
  \end{center}
\end{figure}

\noindent {\bf Remarque. } Il y a une description analogue
en terme d'\'eventails
des \'eclatements le long de centres toriques de dimension
sup\'erieure \`a un (\cite{Oda88} page 38).

\medskip

\noindent {\bf Convention. } La trace sur
la sph\`ere unit\'e $S^{n-1}$ d'un \'eventail $\Delta$ 
de $\RR ^n$ d\'efinit un d\'ecoupage de
$S^{n-1}$. Toutes les figures de ce texte repr\'esentent la trace
sur la sph\`ere $S^{n-1}$ d'une partie d'un \'eventail $n$-dimensionnel.
Si la figure repr\'esente une application torique
entre deux vari\'et\'es toriques, les parties non repr\'esent\'ees
des \'eventails ne sont pas 
modifi\'ees par l'application.

\subsection{Th\'eorie de Mori torique, d'apr\`es Miles Reid}

Si $X$ est une vari\'et\'e projective, on note 
$$\Nun_1(X) = \{ \sum_i a_i C_i \, | \, a_i \in \QQ \, , \, 
C_i \, \, \mbox{courbe irr\'eductible de}\,  \, X \}/\equiv $$
o\`u $\equiv$ d\'esigne l'\'equivalence num\'erique. 
Le c\^one de Mori, ou c\^one des courbes effectives,
est le sous-c\^one de $\Nun_1(X)$ d\'efini par
$$ \NE (X) = \{ Z\in \Nun_1(X) \, | \, Z \equiv 
\sum_i a_i C_i \, , \, a_i\geq 0\}.$$
Si $X$ est torique, alors (voir \cite{Rei83}) $\NE (X)$ est poly\'edral, 
engendr\'e par les classes des
courbes invariantes. De plus, pour
toute ar\^ete $R$ de $\NE (X)$, il y a
une contraction $\varphi_R~: X \to Y$ de $X$ sur une vari\'et\'e torique
projective \'eventuellement singuli\`ere $Y$ telle que 
les courbes irr\'eductibles contract\'ees par $\varphi_R$ sont exactement
celles dont la classe dans $\Nun_1(X)$ appartient \`a
$R$.

\medskip

\noindent {\bf Convention :} pour
une ar\^ete $R$ de $\NE (X)$, nous dirons 
que la contraction extr\'emale $\varphi_R$ est {\em Mori extr\'emale}
si et seulement si $-K_X\cdot R > 0$.  

\subsubsection{Comportement du c\^one de Mori.}
Soit $\varphi_R~: X \to Y$ une contraction extr\'emale
et $$(\varphi_R)_*~: \Nun_1(X) \to \Nun_1(Y) = \Nun_1(X)/(R+(-R))$$
la projection induite par $\varphi_R$.
Alors, le c\^one de Mori $\NE (Y)$
est \'egal \`a $(\varphi_R)_*(\NE (X))$.
En particulier, les ar\^etes de $\NE (Y)$ sont images
d'ar\^etes de $\NE (X)$.

\subsubsection{Description des contractions extr\'emales.}  
Soit $\omega$ une courbe invariante extr\'emale 
(i.e. $R:=\QQ^+[\omega]$ est une ar\^ete de $\NE (X)$)
que l'on peut supposer, quitte \`a renum\'eroter 
les $e_i$, d\'efinie par le ``mur'' 
$<e_1,\ldots,e_{n-1}>$, face 
de deux c\^ones maximaux $$<e_1,\ldots,e_{n-1},e_n>\, \mbox{et} \,
<e_1,\ldots,e_{n-1},e_{n+1}>.$$ 
Comme $X$ est lisse,
il y a une unique relation $$e_{n+1} + e_n + \sum_{i=1}^{n-1}a_i e_i =0$$
o\`u les $a_i$ sont des entiers. 
Rappelons ici que le fibr\'e normal de $\omega$
dans $X$ est isomorphe \`a $\displaystyle{ \bigoplus _{i=1}^{n-1} 
{\mathcal O}_{\PP ^1}(a_i)}$, ceci \'etant vrai m\^eme si
$\omega$ n'est pas extr\'emale.

On note
$$\alpha = \card \{ i \in [1,\ldots,n-1] \, |\, a_i <0\}
\,\,  \mbox{et} \, \, 
\beta = \card \{ i \in [1,\ldots,n-1] \, |\, a_i \leq 0\}.$$ 

Avec ces notations, $\varphi_R$ est birationnelle 
si et seulement si $\alpha \neq 0$, le lieu exceptionnel $A(R)$
de $\varphi_R$ dans $X$ est alors de dimension $n-\alpha$,
son image $B(R)=\varphi_R (A(R))$ dans $Y$ est de 
dimension $\beta-\alpha$,
et la restriction de $\varphi_R$ \`a $A(R)$ est un morphisme plat \`a fibres
des espaces projectifs \`a poids de dimension $n-\beta$.
Si $\alpha =0$, $Y$ est projective lisse de
dimension $\beta$ et $\varphi_R$ est une fibration lisse.
(Attention~: les coquilles de l'\'enonc\'e (2.5) dans \cite{Rei83} sont
corrig\'ees dans \cite{Oda88}, page 111).

\subsection{Un crit\`ere de projectivit\'e torique}

Nous donnons ici une premi\`ere application de la th\'eorie de 
Mori torique.

\begin{prop} Soit $X$ une vari\'et\'e torique, $Z$ une sous-vari\'et\'e
torique de $X$. On suppose que $B_Z(X)$ est
projective.

Alors $X$ est projective si et seulement si
la classe dans $\NE (B_Z(X))$ d'une courbe contenue
dans une fibre de l'\'eclatement $\pi~:B_Z(X) \to X$
est extr\'emale.  

\end{prop} 

\noindent {\bf Remarque :} comme les fibres non triviales de
$\pi$ sont des espaces projectifs (lisses), 
l'ar\^ete engendr\'ee par la classe dans $\NE (B_Y(X))$ d'une courbe contenue
dans une fibre de l'\'eclatement $\pi~:B_Y(X) \to X$
est ind\'ependante des choix de la courbe   
et de la fibre.
 
\medskip 

\noindent {\bf D\'emonstration de la proposition.} Le sens direct 
est un fait g\'en\'eral bien connu~:
soit $F$ une courbe rationnelle contenue
dans une fibre non triviale de $\pi$. Si $X$ est projective, soit
$H$ un diviseur ample sur $X$. Si $F\equiv C_1+C_2$ dans 
$\NE (B_Y(X))$, alors 
$$0=\pi^*H\cdot F = \pi^*H \cdot C_1+\pi^*H \cdot C_2,$$
donc $$\pi^*H\cdot C_1=\pi^*H\cdot C_2 =0.$$
Ainsi les $C_i$
sont contenues dans une fibre de $\pi$, donc
num\'eriquement proportionnelles \`a $F$, autrement dit $F$ est 
extr\'emale.
R\'eciproquement, si $F$ est extr\'emale, la contraction 
extr\'emale $\varphi$ associ\'ee est birationnelle de lieu exceptionnel
\'egal au diviseur exceptionnel de $\pi$~:
en effet, la classe de toute courbe incluse dans une fibre de $\pi$
appartient \`a l'ar\^ete $\QQ^+ [F]$, donc toutes les fibres
non triviales de $\pi$ sont contract\'ees par $\varphi$,
et toute courbe contract\'ee par $\varphi$ est incluse
dans $E$ car $E\cdot F < 0$. Ainsi, 
toute fibre de $\pi$
est contenue dans une fibre de $\varphi$, et les fibres 
de $\varphi$ \'etant des espaces projectifs \`a poids, 
la restriction de $\pi$ \`a toute fibre de $\varphi$
la contracte sur un point, c'est donc que
$\varphi=\pi$ et 
$X$ est projective.\finpreuve

\medskip

\noindent {\bf Remarque :} cet \'enonc\'e 
est un \'enonc\'e torique,
faux pour des vari\'et\'es projectives quelconques 
(voir \cite{Bon96} {\S}2.4) !! 

\section{Un crit\`ere de rigidit\'e}

Le r\'esultat suivant n'utilise que partiellement
l'hypoth\`ese $X$ torique, il serait int\'eressant 
de conna\^\i tre le degr\'e de g\'en\'eralit\'e de cet \'enonc\'e.

\begin{prop}
Soit $X$ une vari\'et\'e torique, $C$ une courbe
torique de $X$. On suppose que $B_C(X)$ est
projective et que $C$ se d\'eforme dans $X$.
Alors $X$ est projective.
\end{prop}

\noindent {\bf Remarque :} dans cet \'enonc\'e, les d\'eformations de 
$C$ sont n\'ecessairement non toriques. 

\medskip 

\noindent {\bf D\'emonstration de la proposition.} 
Comme $X$ est torique, un diviseur sur $X$ est ample
si et seulement si son intersection avec toute
courbe effective est strictement positive.
De l\`a, si $H$ est un diviseur ample sur $B_C(X)$,
et si $\pi~: B_C(X) \to X$ est l'\'eclatement, le
diviseur $H' := \pi (H)$ est ample sur $X$~: en effet,
si $\omega$ est une courbe irr\'eductible de $X$ distincte 
de $C$, alors $H'\cdot \omega \geq H\cdot \tilde{\omega} >0$
o\`u $\tilde{\omega}$ est la transform\'ee stricte de $\omega$.
Comme $C$ se d\'eforme en une courbe $\omega _C$,
on a $H'\cdot C=H'\cdot \omega_C >0$ d'apr\`es ce qui pr\'ec\`ede.
Ainsi, $H'$ est ample sur $X$.\finpreuve

\medskip 

\noindent {\bf Remarque :} la d\'emonstration pr\'ec\'edente
et le crit\`ere de Kleiman \cite{Kle66} montrent
que cette proposition est vraie si on suppose
que $X$ est une vari\'et\'e analytique compacte 
ne poss\'edant pas de courbes irr\'eductibles homologues
\`a $0$ et telle que $\NE (X)$ soit ferm\'e dans $\Nun_1 (X)$.
 
\medskip

La proposition~2 poss\`ede le corollaire imm\'ediat
suivant~:

\begin{cor} Sous les hypoth\`eses de la proposition~2 et
si $\chi (N_{C/X})= -K_X \cdot C+n-3 >0$, alors $X$ est projective.
\end{cor} 

Signalons pour conclure ce paragraphe que ce crit\`ere
de rigidit\'e poss\`ede une application int\'eressante en vue
d'un th\'eor\`eme de Moishezon torique effectif en dimension
trois - voir \cite{Bon99}.

\section{Trois exemples de contractions extr\'emales}

Dans cette partie, la situation qui nous int\'eresse
est la suivante~: $X$ est une
vari\'et\'e torique non projective 
contenant
une courbe $C$ torique de sorte que 
$\tilde X~:=B_C(X)$ soit projective.
On note $\pi~: \tilde X \to X$
l'\'eclatement de $X$ le long de $C$, $E$ le diviseur exceptionnel
de $\pi$.

\medskip

\noindent {\bf Question~:} que peut-on dire des contractions 
{\em Mori extr\'emales} 
sur $\tilde X$ dont le lieu exceptionnel rencontre $E$ ?

\medskip

Nous donnons trois exemples diff\'erents de 
situations possibles.

\subsection{Exemple~1~: ``r\'eduction triviale''.}

Supposons donn\'ee une vari\'et\'e torique non projective $X'$,
contenant une courbe $C'$ torique de sorte que 
$Y:=B_{C'}(X')$ soit projective~; on note $\pi'$ l'\'eclatement 
$Y \to X'$.
Soit alors $X$ la vari\'et\'e torique obtenue 
en \'eclatant $X'$ en un point fixe $p$ 
appartenant \`a $C'$.
Comme $X'$, la vari\'et\'e $X$ est non projective, et la vari\'et\'e
$\tilde X := B_C(X)$
o\`u
la courbe $C$ est 
la transform\'ee stricte de $C'$
est projective~: en effet,
$\tilde X$
est la vari\'et\'e obtenue en \'eclatant $Y$ le long de
la sous-vari\'et\'e 
de codimension deux $(\pi') ^{-1}(p)$.
D'apr\`es la proposition~1, l'\'eclatement
$\tilde X \to Y$ est une contraction Mori extr\'emale
$\varphi _{\QQ ^+ [\omega]}$ dont le lieu exceptionnel 
rencontre $E$.

Autrement dit, on a un diagramme 
commutatif~:

\centerline{
\xymatrix{ \tilde{X}=B_C(X) \ar[d]_{\pi} \ar[r]^{\varphi _{\QQ ^+ [\omega]}}
  & Y = B_{C'}(X') \ar[d]^{\pi '} \\
X = B_p(X') \ar[r]^{\varphi} & X'
}
}

La figure~2 est la version ``g\'eom\'etrique''
de ce qui pr\'ec\`ede.

\medskip

\noindent {\bf Convention. } Si le couple $(X,C)$ 
est obtenu par le proc\'ed\'e d\'ecrit pr\'ec\'edemment,
nous dirons que ce couple poss\`ede une {\em r\'eduction 
triviale}.

\bigskip

%cette figure est exportee en 50,20
\begin{figure}[hbtp]
  \begin{center}
    \leavevmode
    \input{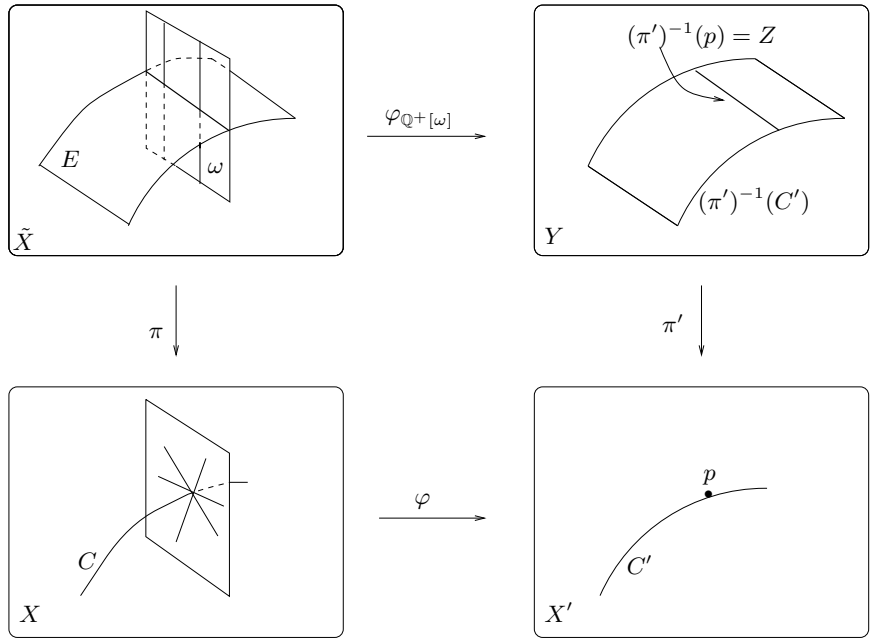}
    \caption{R\'eduction triviale}
    %\label{fig:}
  \end{center}
\end{figure}

\bigskip

\subsection{Exemple~2~: ``flip interdit''.}

Nous d\'emontrerons au {\S}4.2 la proposition suivante~:

\begin{prop}
Soit $Y$ une vari\'et\'e torique projective de dimension $n$
contenant une sous-vari\'et\'e torique $Z$ 
v\'erifiant 
$$(Z, N_{Z/Y})\simeq (\PP ^{n-2}, 
{\mathcal O}_{\PP ^{n-2}}(-1)^{\oplus 2}),$$
et soit $X$ la vari\'et\'e torique obtenue 
en \'eclatant $Y$ le long de $Z$
puis en contractant les $\{*\} \times \PP ^{n-2}$
du diviseur exceptionnel 
$E\simeq \PP^1 \times \PP ^{n-2}$~:

\medskip

\centerline{
\xymatrix{  &\tilde{X}=B_Z(Y) \ar[dl]_\pi 
\ar[dr]^{\varphi _{\QQ ^+ [\omega]}} \\
X & & Y
}
}

Alors $X$ est projective si et seulement si
la classe d'une courbe 
incluse dans $Z$ est  
ex\-tr\^e\-ma\-le dans le c\^one de
Mori des courbes effectives de $Y$.
\end{prop}

En particulier, si la classe d'une courbe 
incluse dans $Z$ n'est pas  
ex\-tr\^e\-ma\-le dans le c\^one de
Mori des courbes effectives de $Y$, la vari\'et\'e $X$ construite
est non projective, le devient apr\`es 
\'eclatement le long d'une courbe $C$
v\'erifiant $$ (C, N_{C/X}) \simeq (\PP ^1, 
{\mathcal O}_{\PP^1}(-1)^{\oplus n-1}),$$
et $\varphi _{\QQ ^+ [\omega]}$ est une contraction Mori extr\'emale
dont le lieu exceptionnel 
est \'egal \`a celui de l'\'eclatement $\pi~: \tilde X \to X$.

\medskip

\noindent {\bf Convention. } Si une vari\'et\'e $X$ non projective est 
obtenue par le proc\'ed\'e d\'ecrit pr\'ec\'edem\-ment,
nous dirons que cette vari\'et\'e est obtenue {\em via}
un {\em flip interdit} d'une vari\'et\'e projective~: 
\`a partir de la dimension quatre, cette transformation 
de $Y$ vers $X$ est un flip
si la classe d'une courbe 
incluse dans $Z$ est  
ex\-tr\^e\-ma\-le dans le c\^one de
Mori des courbes effectives de $Y$.

\subsection{Exemple~3~: ``transformation \'el\'ementaire''.}

Soit $Y$ une vari\'et\'e torique projective de dimension $n$
contenant un diviseur torique $\bar E$. Supposons
que $\bar E$ soit isomorphe \`a $\PP ({\mathcal E})$
o\`u ${\mathcal E}$ est un fibr\'e vectoriel de rang $n-1$
sur $\PP ^1$ et notons $\lambda~: \bar E= \PP ({\mathcal E}) \to \PP ^1$
la projection. Supposons que le fibr\'e normal
de $\bar E$ dans $Y$ soit de la forme $\lambda ^* {\mathcal O}_{\PP ^1}(a)$
(autrement dit, il est trivial en restriction aux fibres de $\lambda$)
et soit ${\mathcal E'}$ un sous-fibr\'e de ${\mathcal E}$
de rang $n-2$ de sorte que la sous-vari\'et\'e $Z = \PP ({\mathcal E'})$
soit torique. Soit $\tilde{X} = B_Z(Y)$ la vari\'et\'e torique
obtenue en \'eclatant $Y$ le long de $Z$. 
D'apr\`es la proposition~1, l'\'eclatement
$\tilde X \to Y$ est une contraction Mori extr\'emale
$\varphi _{\QQ ^+ [\omega]}$.
De plus, la transform\'ee stricte $E$ de $\bar E$
est isomorphe \`a $\bar E$ (en particulier, 
$\lambda$ induit une fibration $\tilde{\lambda}~: E \to \PP ^1$) 
et son fibr\'e normal
$N_{E / \tilde X}$ est \'egal \`a ${\mathcal O}_{E}(-1)\otimes
\tilde{\lambda}^* {\mathcal O}_{\PP ^1}(a)$.
Il existe donc une vari\'et\'e $X$ torique contenant
une courbe torique $C \simeq \PP ^1$ de sorte que $\tilde X = B_C(X)$ et 
que le diviseur exceptionnel de l'\'eclatement $\pi~: \tilde X \to \PP ^1$
soit \'egal \`a $E$ avec $\pi_{|E} = \tilde{\lambda}$.
Par construction, 
$\varphi _{\QQ ^+ [\omega]}$ est une contraction Mori extr\'emale
dont le lieu exceptionnel $F$
rencontre $E$.
 
Autrement dit, on a un diagramme 
commutatif~:

\centerline{\xymatrix{  &E \subset \tilde{X} \supset F
 \ar[dl]_\pi 
\ar[dr]^{\varphi _{\QQ ^+ [\omega]}} \\
C \subset X & & Y \supset \bar E \supset Z
}
}

\medskip

La figure~3 est la version ``g\'eom\'etrique''
de ce qui pr\'ec\`ede.
 
\medskip

\noindent {\bf Convention. } Si une vari\'et\'e $X$ est 
obtenue par le proc\'ed\'e d\'ecrit pr\'ec\'edem\-ment,
nous dirons que cette vari\'et\'e est obtenue
{\em via} une {\em transformation \'el\'ementaire}
d'une vari\'et\'e projective.

\medskip

%cette figure est exportee en 50,20
\begin{figure}[hbtp]
  \begin{center}
    \leavevmode
    \input{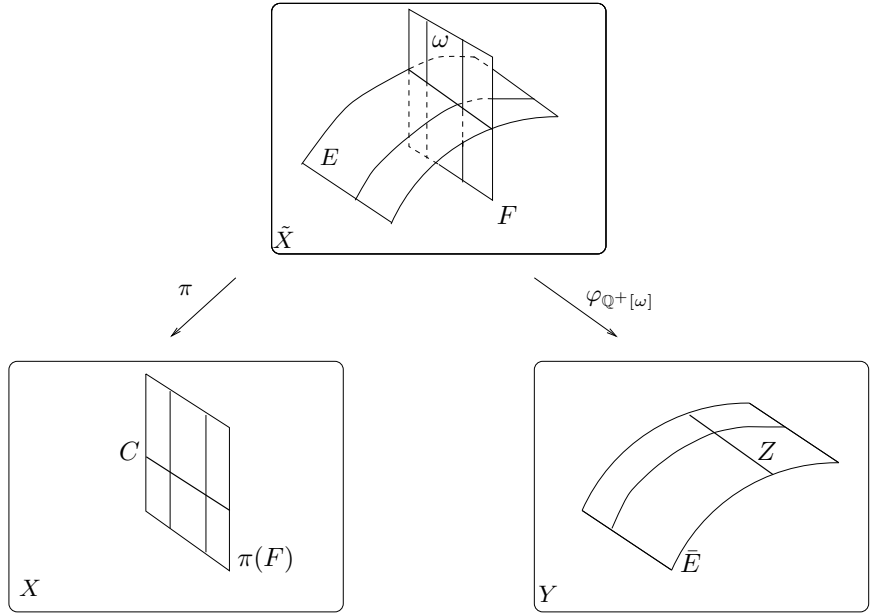}
    \caption{Transformation \'el\'ementaire.}
    %\label{fig:}
  \end{center}
\end{figure} 

\subsubsection*{Lien avec les transformations \'el\'ementaires
de Maruyama.}
Remarquons que dans ce qui pr\'ec\`ede, l'image dans $X$ du diviseur
$F$ est encore un diviseur de la forme 
$\PP ({\mathcal F})$,
o\`u ${\mathcal F}$ est un fibr\'e vectoriel de rang $n-1$
sur $\PP ^1$, dont le fibr\'e normal
dans $X$ est encore trivial en restriction aux fibres de la
projection sur $\PP ^1$. Cette construction
est un cas particulier de transformation \'el\'ementaire, 
application birationnelle
utilis\'ee par Maruyama \cite{Mar82}
pour construire des familles de
fibr\'es vectoriels 
sur une vari\'et\'e analytique ou alg\'ebrique (en dimension deux, 
cette application birationnelle joue un r\^ole
essentiel dans la
classification des surfaces rationnelles).

Plus pr\'ecis\'ement (voir \cite{Mar82} {\S}1),
si ${\mathcal E}$ est un fibr\'e vectoriel
sur une vari\'et\'e $S$, si $T$ est un diviseur (lisse) de $S$ 
et si ${\mathcal E'}$ est un sous-fibr\'e vectoriel strict
de ${\mathcal E}_{|T}$, une transformation \'el\'ementaire 
consiste \`a \'eclater la vari\'et\'e $\PP ({\mathcal E})$
le long de la sous-vari\'et\'e (de codimension au moins deux)
$\PP ({\mathcal E'})$, et \`a remarquer que la transform\'ee stricte 
du diviseur $\PP ({\mathcal E}_{|T})$ est le diviseur exceptionnel
d'un \'eclatement sur une vari\'et\'e de la forme $\PP ({\mathcal F})$
o\`u ${\mathcal F}$ est encore un fibr\'e vectoriel sur $S$, de
m\^eme rang que ${\mathcal E}$.

Dans notre situation, notons $S$ l'espace total 
du fibr\'e ${\mathcal O}_{\PP ^1}(a)$, 
$q~: S \to \PP^1$ la projection et $U$ le produit fibr\'e
(torique) de $\lambda$ et $q$ au dessus de $\PP^1$~:  

\centerline{
\xymatrix{ U \ar[d]  \ar[r]
  & S \ar[d]^{q}  \\
\bar E \ar[r]^{\lambda} & \PP ^1
}
}
L'hypoth\`ese faite sur $N_{\bar E/Y}$
implique qu'il y a un voisinage torique $V$ de $\bar E$ dans $Y$
isomorphe \`a $U = \PP (q^*{\mathcal E})$. 
Enfin, si $T$ d\'esigne la section nulle de $q$ dans $S$ et si  
$\tilde V = 
(\varphi _{\QQ ^+ [\omega]})^{-1}(V)$, alors le diagramme

\centerline{\xymatrix{  &\tilde{V} \ar[dl]_{\pi _{| \tilde V}} 
\ar[dr]^{(\varphi _{\QQ ^+ [\omega]})_{| \tilde V} } \\
\pi (V) & & V
}
}

est une transformation \'el\'ementaire de Maruyama. 

A la lumi\`ere de ce qui pr\'ec\`ede, nous pouvons dire que 
l'exemple~3 est une ``compactification \'equivariante 
d'une transformation \'el\'ementaire torique''.

\section{Classification de certaines contractions Mori extr\'emales}

Dans toute cette partie, $X$ est une
vari\'et\'e torique non projective de dimension $n$, d'\'eventail 
$\Delta$
contenant
une courbe $C$ torique de sorte que 
$\tilde X~:=B_C(X)$ soit projective, d'\'eventail 
$\tilde{\Delta}$. On note $\pi~: \tilde X \to X$
l'\'eclatement de $X$ le long de $C$, $E$ le diviseur exceptionnel
de $\pi$.
Quitte \`a renum\'eroter 
les faces de dimension un de $\Delta$, 
on suppose que $C$ est donn\'ee par un c\^one $<e_1,\ldots,e_{n-1}>$ 
de
dimension $n-1$, 
face de deux c\^ones maximaux $<e_1,\ldots,e_{n-1},e_n>$
et $<e_1,\ldots,e_{n-1},e_{n+1}>$.

Le but de cette partie est de classifier les contractions 
{\em Mori extr\'emales} 
de $\tilde X$ dont le lieu exceptionnel rencontre $E$.
Nous y d\'emontrons le r\'esultat suivant~:

\bigskip

\noindent {\bf Th\'eor\`eme~A. } {\em Soit $X$ une
vari\'et\'e torique non projective
contenant une courbe $C$ torique de sorte que 
$\tilde X~:=B_C(X)$ soit projective, soit $\omega$
dans $\tilde X$ une courbe Mori extr\'emale 
rencontrant le diviseur exceptionnel $E$ et
$\varphi_{\QQ ^+[\omega]}~: \tilde X \to Y$
la contraction extr\'emale associ\'ee.
Alors~:
\begin{enumerate}

\item[$\bullet$] soit le couple $(X,C)$ poss\`ede une r\'eduction 
triviale,

\item[$\bullet$] soit $X$ est obtenue {\em via}
un flip interdit d'une vari\'et\'e torique projective,

\item[$\bullet$] soit $X$ est obtenue
{\em via} une transformation \'el\'ementaire
d'une vari\'et\'e torique projective.

\end{enumerate}

Dans tous les cas, la contraction extr\'emale
$\varphi_{\QQ ^+[\omega]}$ est de la forme d\'ecrite
dans les exemples 1, 2 et 3 respectivement.  

}

\bigskip

\noindent {\bf Commentaires.} 
\begin{enumerate}
\item[(i)] Dans tous les cas, l'image $Y$ de la
contraction extr\'emale $\varphi _{\QQ ^+ [\omega]}$
est lisse et $\varphi _{\QQ ^+ [\omega]}$ est 
l'\'eclatement de $Y$ le long d'un centre lisse~; 
ce ph\'enom\`ene est remarquable.
\item[(ii)] 
Dans le deuxi\`eme cas, $E\cdot \omega =-1$,
dans les deux autres $E\cdot \omega =1$.
\item[(iii)] Tous ces exemples se produisent en
toutes dimensions sup\'erieures \`a trois,
nous donnerons des constructions explicites dans la partie~5
avec ``petits'' groupes de Picard. 
\item[(iv)] Cet \'enonc\'e n'est pertinent que s'il
existe une courbe Mori extr\'emale rencontrant $E$,
hypoth\`ese que nous ne savons pas v\'erifier en g\'en\'eral.
Nous donnons dans la proposition suivante deux
conditions sous lesquelles cette hypoth\`ese est v\'erifi\'ee.

\end{enumerate}

\begin{prop}
Soit $X$ une vari\'et\'e torique non projective
contenant une courbe $C$ torique de sorte que 
$\tilde X~:=B_C(X)$ soit projective, soit $E$
le diviseur exceptionnel de l'\'eclatement et
$\rho (X)$ le nombre de Picard de $X$.
Alors $\tilde X$ poss\`ede une courbe Mori extr\'emale 
rencontrant $E$ d\`es que $\tilde X$ v\'erifie
l'une des conditions suivantes~:
\begin{enumerate}
\item[(i)] la vari\'et\'e $\tilde X$ est de Fano,

\item[(ii)] la vari\'et\'e $\tilde X$ poss\`ede au moins $\rho (X)$
contractions Mori extr\'emales ({\em i.e} $\NE (\tilde X)$
poss\`ede au moins $\rho (X)$ ar\^etes $R$ v\'erifiant
$-K_{\tilde X}\cdot R > 0$).

\end{enumerate}
\end{prop}

\noindent {\bf D\'emonstration de la proposition.}

Pour (i)~: si $\tilde X$ est de Fano, soit $C$ 
une courbe contenue dans une fibre de l'\'eclatement.
Alors, 
$$ [C] \equiv \sum_{i=1}^p a_i [C_i] $$
o\`u les $C_i$ sont des courbes Mori extr\'emales
et les $a_i$ des rationnels positifs.
Comme $E \cdot C <0$, l'une au moins des 
$C_i$ v\'erifie $E \cdot C_i <0$.

Pour (ii)~: si par l'absurde, il n'y a pas de courbe Mori extr\'emale 
rencontrant $E$, alors l'hyperplan 
$$ E_{=0} = \{Z \in \Nun _1(\tilde X) \, | \,
E\cdot Z =0 \} $$
contient $\rho (X)$ ar\^etes de $\NE (\tilde X)$. 
Or $\NE (\tilde X)$ est un c\^one convexe d'int\'erieur non vide 
dans l'espace vectoriel $\Nun _1(\tilde X)$ de dimension 
$\rho (X)+1$. Il s'en suit que $ E_{=0}$
est un hyperplan d'appui de $\NE (\tilde X)$. Ceci est absurde 
car il y a dans $\NE (\tilde X)$ des courbes d'intersection
strictement po\-sitive avec $E$ (toute courbe transverse \`a $E$)
et des courbes d'intersection strictement n\'egative avec $E$
(toute courbe contenue dans une fibre de l'\'eclatement).\finpreuve

\newpage

La preuve pr\'ec\'edente (pour (i)) et le th\'eor\`eme~A
implique de suite le r\'esultat suivant.

\medskip

\noindent {\bf Corollaire~B. } {\em Soit $X$ une
vari\'et\'e torique non projective
contenant une courbe $C$ torique de sorte que 
$\tilde X~:=B_C(X)$ soit une vari\'et\'e de Fano.
Alors $X$ est obtenue {\em via}
un flip interdit d'une vari\'et\'e torique projective.
}

\medskip

En dimension $3$, on en d\'eduit~:

\medskip

\noindent {\bf Corollaire~C. } {\em 
Soit $X$ une
vari\'et\'e torique de dimension $3$
contenant une courbe $C$ torique de sorte que 
$\tilde X~:=B_C(X)$ soit une vari\'et\'e de Fano.
Alors $X$ est projective.
}

\medskip

\noindent {\bf D\'emonstration.} Si $X$ n'est pas projective,
le corollaire~B affirme que $X$ est obtenue {\em via}
un flip interdit d'une vari\'et\'e torique projective.
Comme $-K_{\tilde X}$ est ample, $-K_X$ est d'intersection
strictement positive avec toutes les courbes de $X$
distinctes de $C$. Or, comme 
$$N_{C/X} = {\mathcal O}_{\PP ^1}(-1)^{\oplus 2},$$
on a $-K_X\cdot C =0$.
Mais alors, si $A$ est un diviseur nef dans une vari\'et\'e 
torique, d'intersection strictement positive sur
toutes les courbes toriques \`a l'exception d'une seule not\'ee $C$,
et si $D$ est un diviseur lisse v\'erifiant $D\cdot C =1$,
alors $A + \varepsilon D$ est un $\QQ $-diviseur 
strictement positif sur toutes les courbes toriques, 
donc ample,
d\`es 
que $\varepsilon$ est un rationnel strictement positif
assez petit. Ainsi, $X$ est projective.\finpreuve   

\medskip

\noindent {\bf Remarque~:} le corollaire~C d\'ecoule aussi de
la classification des vari\'et\'es toriques de Fano de dimension 
$3$ \cite{Bat82} \cite{WWa82}. En effet, cette derni\`ere 
montre que $\rho (\tilde X) \leq 5$, donc $\rho (X) \leq 4$.
Or il n'y a qu'une seule vari\'et\'e torique 
non projective de dimension $3$ dont le rang du groupe de Picard est 
inf\'erieur ou \'egal \`a $4$ (celle dont l'\'eventail est rappel\'e
au {\S}5.1.1 de ce texte)
\cite{Oda88} et il est facile de voir que cette derni\`ere
ne devient pas de Fano apr\`es \'eclatement le long d'une
courbe. La preuve du corollaire~C \'evite d'utiliser cette classification.

\bigskip

{\bf Le reste de cette partie est consacr\'e
\`a la d\'emonstration du th\'eor\`eme~A.
}

\medskip

Supposons donc l'existence d'une 
courbe torique $\omega$
de $\tilde X :=B_C(X)$ Mori extr\'emale 
v\'erifiant $E\cap \omega \neq \emptyset$.
L'analyse distingue les cas $E\cdot\omega \geq 0$
et $E\cdot\omega < 0$. 

\subsection{Le cas o\`u $E\cdot\omega \geq 0$.}

\subsubsection{Quelques r\'eductions}

\begin{lemm}
Sous les hypoth\`eses pr\'ec\'edentes, la courbe
$\omega$ n'est pas incluse dans $E$.
\end{lemm}

\noindent {\bf D\'emonstration.}
Supposons au contraire que $\omega \subset E$.
Comme $E\cdot\omega \geq 0$, $\omega$ n'est pas incluse 
dans une fibre de
$\pi$, donc $\pi_{|\omega} : \omega \to C$
est un isomorphisme d'apr\`es 1.2.
De l\`a, 
$$ K_X \cdot C= \pi^*K_X\cdot \omega = 
(K_{\tilde X}-(n-1)E) \cdot \omega$$
et donc 
$$\chi (N_{C/X})= -K_X \cdot C+n-3 
= (-K_{\tilde X}+(n-1)E) \cdot \omega +n-3 >0.$$
D'apr\`es le corollaire~1, ceci implique que $X$ est projective, ce qui est 
exclu par hypoth\`ese.\finpreuve

\medskip

Ainsi, $\omega$ est ``transverse'' \`a $E$ et 
quitte \`a permuter les 
$e_1,\ldots,e_{n-1}$, on peut 
supposer que $\omega =<e_1,\ldots,e_{n-2},e_n>$
(et donc $E\cdot \omega =1$ ).
Le mur 
$<e_1,\ldots,e_{n-2},e_n>$ d\'efinissant $\omega$
est face de deux c\^ones maximaux $$<e_1,\ldots,e_{n-2},e_n,e>\,\, 
\mbox{ et }
<e_1,\ldots,e_{n-2},e_n,e'>$$ pour un certain $e'$,
g\'en\'erateur entier d'une face de dimension un
de $\tilde{\Delta}$. Remarquons que $e'$ peut \^etre
\'egal \`a $e_{n+1}$ et est distinct de $e_{n-1}$ - voir
la figure~4.

%cette figure est exportee a 60, 16
\begin{figure}[hbtp]
  \begin{center}
    \leavevmode
    \input{fig3.pstex_t}
    \caption{}
    %\label{fig:}
  \end{center}
\end{figure}

Il y a alors une relation 
$$ e' + e + \sum_{i=1} ^{n-2} b_i e_i + b_n e_n = 0,$$
o\`u les $b_i$ et $b_n$ sont des entiers.

\newpage

\begin{lemm} Les entiers $b_i$ sont n\'egatifs ou nuls
pour tout $i=1,\ldots,n-2,n$.
\end{lemm}  

\noindent {\bf D\'emonstration.}
Sinon, on a $b_i >0$ pour un certain $i$, et d'apr\`es Reid 
\cite{Rei83} Cor. 2.10(i),
les deux courbes 
$$ \omega':= <e_1,\ldots,\hat{e_i},\ldots,e_{n-2},e_n,e>\,\, \mbox{ et }
\,\, \omega'':= 
<e_1,\ldots,\hat{e_i},\ldots,e_{n-2},e_n,e'>$$ 
(respectivement $<e_1,\ldots,e_{n-2},e>$ et
$<e_1,\ldots,e_{n-2},e'>)$ si $i\leq n-2$
(respectivement si $i=n$)
sont sur l'ar\^ete
$\QQ ^+ [\omega]$.

Mais alors $\omega':= <e_1,\ldots,\hat{e_i},\ldots,e_{n-2},e_n,e>$
(respectivement $<e_1,\ldots,e_{n-2},e>$)
est une courbe Mori extr\'emale incluse dans $E$
et d'intersection positive avec $E$, ce qui n'est pas possible
d'apr\`es le lemme pr\'ec\'edent.
\finpreuve

\medskip

Comme la quantit\'e $$-K_{\tilde X}\cdot \omega  = 2 + \sum_{i=1} ^{n-2} b_i
+ b_n$$ est strictement positive, les deux lemmes pr\'ec\'edents
et la description des contractions extr\'e\-ma\-les rappel\'ees en 1.3 
impliquent que seules deux possibilit\'es ont lieu {\em a priori}.

\begin{enumerate}

\item [a)] Tous les entiers $b_i$, $1\leq i \leq n-2$
et $b_n$ sont nuls et dans ce cas
la contraction extr\'emale $\varphi_{\QQ ^+[\omega]}$
est une fibration $\varphi_{\QQ ^+[\omega]}~: \tilde X \to Y$ 
de fibres isomorphes \`a $\PP ^1$.   

\item [b)] Un et un seul des entiers $b_i$ et $b_n$ vaut $-1$,
les autres \'etant nuls et dans ce cas
la contraction extr\'emale $\varphi_{\QQ ^+[\omega]}$ est birationnelle
divisorielle, ses fibres non triviales \'etant isomorphes \`a $\PP ^1$.

\end{enumerate}

Montrons que le cas  a) est impossible~: si
$\varphi_{\QQ ^+[\omega]}~: \tilde X \to Y$ est une fibration,
ses fibres ont une intersection \'egale \`a $1$ avec le diviseur
$E$ et donc $E$ est isomorphe \`a $Y$. On en d\'eduit
que le groupe de Picard de $Y$ est $\ZZ ^2$, donc celui de 
$\tilde X$ est $\ZZ ^3$, donc celui de $X$ est
$\ZZ ^2$. Mais alors $X$ est projective d'apr\`es la classification
\cite{Kle88}.

\subsubsection{Analyse du cas b).}
Deux situations diff\'erentes peuvent se produire~: 

\begin{enumerate}
\item[(i)] $b_i=0$ pour tout $1\leq i \leq n-2$
et $b_n=-1$, on a 
alors $e' + e  - e_n = 0$,
\item[(ii)] l'un des $b_i$, $1\leq i \leq n-2$ est non 
nul, on peut supposer que $b_1=-1$ et on a 
alors $ e' + e  - e_1 = 0$.
\end{enumerate}

Commen\c cons par traiter {\bf la situation (i)}~: nous
montrons que $(X,C)$ poss\`ede une {\em r\'eduction triviale}.

\medskip

%cette figure est exportee en 50,20
\begin{figure}[hbtp]
  \begin{center}
    \leavevmode
    \input{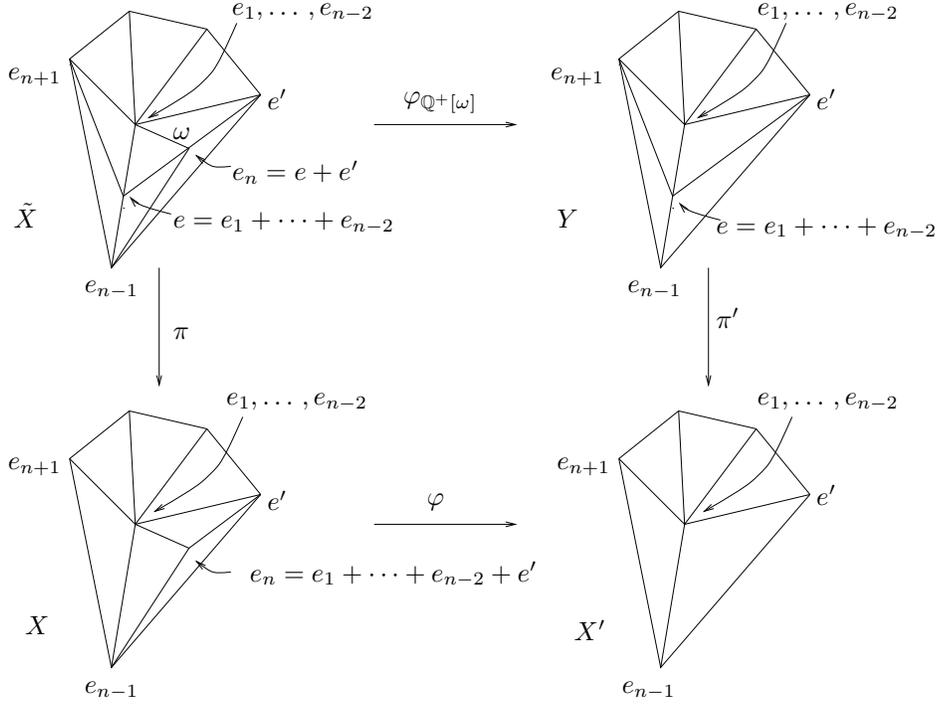}
    \caption{R\'eduction triviale}
    %\label{fig:}
  \end{center}
\end{figure}

\medskip

D'apr\`es Reid \cite{Rei83} Cor. 2.10(ii), 
pour tout $i$ fix\'e, $1\leq i \leq n-2$,
l'\'etoile du c\^one $n-2$ dimensionnel 
$<e_1,\ldots,\hat{e_i}, \ldots,e_{n-2},e_n>$ 
est constitu\'ee des $4$ c\^ones
maximaux $$<e_1,\ldots,e_i, \ldots,e_{n-2},e_n,e> \, , \,
<e_1,\ldots,e_i, \ldots,e_{n-2},e_n,e'>\, ,\, $$
$$ <e_1,\ldots,\hat{e_i}, \ldots,e_{n-2},e_n,e,\nu_i>\, \mbox{et} \,
<e_1,\ldots,\hat{e_i}, \ldots,e_{n-2},e_n,e',\nu_i>$$  
pour un certain $\nu _i$. 

Or, le
mur $<e_1,\ldots,\hat{e_i}, \ldots,e_{n-2},e_n,e>$
est d\'ej\`a face des $2$ c\^ones
maximaux $$<e_1,\ldots,\hat{e_i}, \ldots,e_{n-2},e_n,e,e_{n-1}>\, \mbox{et} \,
<e_1,\ldots,e_i, \ldots,e_{n-2},e_n,e>.$$
C'est donc que $\nu _i = e_{n-1}$
pour tout $1\leq i \leq n-2$.
On en d\'eduit que l'\'etoile de $e_n$
est constitu\'ee des $2n-2$
c\^ones
maximaux suivants~: $$<e_1,\ldots,\hat{e_i},\ldots,e_{n-1},e,e_n> \,\,
\mbox{ et }
<e_1,\ldots,\hat{e_i},\ldots,e_{n-1},e',e_n>, \,
1\leq i \leq n-1.$$ Il r\'esulte du raisonnement pr\'ec\'edent
-~voir la figure~5 - 
que la contraction $\varphi_{\QQ ^+[\omega]}$
envoie $\tilde X$ sur la vari\'et\'e projective 
lisse $Y$
dont l'\'eventail est obtenu en rempla\c cant dans
$\tilde{\Delta}$ l'\'etoile de $e_n$ par la r\'eunion
des c\^ones maximaux $<e_1,\ldots,\hat{e_i},e_{n-1},e',e>$,
$1\leq i \leq n-1$ et leurs faces.

Autrement dit, 
$\varphi_{\QQ ^+[\omega]}~: \tilde X \to Y$
est l'\'eclatement de $Y$ le long de la sous-vari\'et\'e torique
de codimension deux d\'efinie par le c\^one $<e,e'>$.

De m\^eme, la relation 
$$ e_n = e_1+\cdots + e_{n-1}+e'$$
montre qu'il existe une vari\'et\'e $X'$
non projective lisse telle que $X=B_p(X')$ 
o\`u $p$ est un point fixe dans $X'$ d\'efinie
par le c\^one maximal $<e_1,\ldots,e_{n-1},e'>$.
Enfin, $Y=B_{C'}(X')$ o\`u
$C'$ est la courbe de $X'$ d\'efinie
par le mur $<e_1,\ldots,e_{n-1}>$. Ce qui pr\'ec\`ede montre
que $(X,C)$ poss\`ede une r\'eduction triviale.

\medskip

Traitons {\bf la situation (ii)}~: nous montrons que
$X$ est obtenue
{\em via une transformation \'el\'ementaire}
d'une vari\'et\'e torique projective

Dans ce cas, on a $e_1 = e+e'$, soit
encore 
$e'=-e_2-\cdots-e_{n-1}$.
Un raisonnement analogue au pr\'ec\'edent montre
que l'\'etoile de $e_1$ est constitu\'ee des 
$4n-8$ c\^ones maximaux suivants~:
$$<e_1,\ldots, \hat{e_i}, \ldots, e_{n-1}, e, e_n>,\, 
<e_1,\ldots, \hat{e_i}, \ldots, e_{n-1}, e, e_{n+1}>,$$
$$<e_1,\ldots, \hat{e_i}, \ldots, e_{n-1}, e', e_n>\, \mbox{
et} \, <e_1,\ldots, \hat{e_i}, \ldots, e_{n-1}, e', e_{n+1}>$$
pour $2\leq i \leq n-1$.

Il en r\'esulte - voir la figure~6 -  
que la contraction $\varphi_{\QQ ^+[\omega]}$
envoie $\tilde X$ sur la vari\'et\'e projective 
lisse $Y$
dont l'\'eventail est obtenu en rempla\c cant dans
$\tilde{\Delta}$ l'\'etoile de $e_1$ par la r\'eunion
des c\^ones maximaux $$<e_2,\ldots,\hat{e_i},e_{n-1},e',e,e_n>,\, 
<e_2,\ldots,\hat{e_i},e_{n-1},e',e,e_{n+1}>,\,
2\leq i \leq n-1$$ et leurs faces. 

\bigskip

%cette figure est exportee en 50,20
\begin{figure}[hbtp]
  \begin{center}
    \leavevmode
    \input{fig8.pstex_t}
    \caption{Transformation \'el\'ementaire.}
    %\label{fig:}
  \end{center}
\end{figure} 

\bigskip

Autrement dit, 
$\varphi_{\QQ ^+[\omega]}~: \tilde X \to Y$
est l'\'eclatement de $Y$ le long de la sous-vari\'et\'e torique
de codimension deux d\'efinie par le c\^one $<e,e'>$.
De plus, dans l'\'eventail $\Delta_Y$ d\'efinissant $Y$, l'\'etoile de 
$e$ est constitu\'ee des $2n-2$ c\^ones maximaux 
$$<e_2,\ldots,\hat{e_i},e_{n-1},e',e,e_n>,\, 
<e_2,\ldots,\hat{e_i},e_{n-1},e',e,e_{n+1}>,\,
2\leq i \leq n-1,$$
$$<e_2,\ldots,e_{n-1},e,e_n>\,
\mbox{et} \,<e_2,\ldots,e_{n-1},e,e_{n+1}>$$ et leurs faces.
En particulier, la forme lin\'eaire
$\lambda~: N_{\QQ} \to \QQ$
d\'efinie par $\lambda (e)= \lambda (e_i)= 0$ si $2\leq i \leq n-1$ et
$\lambda (e_n)=1$ (et donc $\lambda (e_{n+1})=-1$ !) 
induit un morphisme surjectif de l'\'eventail de l'\'etoile
de $e$ dans $\Delta_Y$ sur l'\'eventail de $\PP^1$
et donc induit une fibration torique lisse du
diviseur torique $E$ d\'efini par le c\^one $<e>$ sur $\PP^1$,
\`a fibre isomorphe \`a $\PP ^{n-2}$. Enfin, 
comme 
$$e'+e_2 + \sum_{i=3}^{n-1}e_i =0,$$ 
la courbe rationnelle
$<e_3,\ldots,e_{n-1},e,e_{n}>$ (incluse dans une fibre de $\lambda$)
a pour fibr\'e normal dans $Y$ le fibr\'e ${\mathcal O}^{\oplus 2} \oplus
{\mathcal O}(1)^{\oplus n-3}$
et le fibr\'e normal $N_{E/Y}$ est 
donc trivial en restriction aux fibres de $\lambda$.
Tout ceci ach\`eve de montrer que $X$  
est obtenue
{\em via} une transformation \'el\'ementaire
de la vari\'et\'e torique projective $Y$.

\bigskip

\noindent{\bf Remarque~:} dans la situation pr\'ec\'edente,
la classe
d'une courbe incluse dans une fibre de $\lambda~: E \to \PP ^1$
n'est pas extr\'emale dans $\NE (Y)$.

\bigskip

\noindent {\bf D\'emonstration de la remarque.} 
D'apr\`es 1.3, comme $\varphi _{\QQ ^+ [\omega]}$
est une contraction extr\'emale, si une courbe incluse 
dans une fibre de $\lambda~: E \to \PP ^1$
est extr\'emale dans $\NE (Y)$, sa transform\'ee stricte 
l'est dans $\tilde X$, et donc $\pi$ est une contraction extr\'emale
d'apr\`es la proposition~1. Or $X$ n'est pas projective ! \finpreuve
 
\newpage 

\noindent{\bf Remarque~:} dans la situation pr\'ec\'edente, on a
un diagramme~:

\medskip

\centerline{
\xymatrix{ & \tilde{X}=B_C(X) \ar[dl]_{\pi} \ar[dr]
^{\varphi _{\QQ ^+ [\omega]}} &\\
X&  & Y 
}
}

\medskip

o\`u les deux fl\`eches $\pi$ et $\varphi _{\QQ ^+ [\omega]}$
sont des \'eclatements le long de centres lisses.
En fait, on peut it\'erer la construction pr\'ec\'edente~:
\'eclatons $Y$ le long de
la courbe de c\^one $<e_2,\ldots, e_{n-1},e>$, 
la vari\'et\'e projective alors obtenue est elle-m\^eme
l'\'eclatement le long d'un centre lisse
de codimension deux dans une autre vari\'et\'e torique lisse,
et on peut recommencer. Ainsi $X$ fait partie
d'une suite infinie de vari\'et\'es toriques lisses, s'obtenant 
successivement \`a l'aide de transformations \'el\'ementaires. 
Nous ne savons pas 
d\'eterminer de fa\c con syst\'ematique le caract\`ere
projectif ou non de ces vari\'et\'es. Nous verrons dans la partie~5
deux exemples de comportements distincts.   

\subsection{Le cas o\`u $E\cdot \omega < 0$.}  
On d\'emontre la proposition suivante~:

\begin{prop} Soit $X$ une
vari\'et\'e torique non projective. Supposons que~:
\begin{enumerate}
\item [(i)] il existe
une courbe $C$ torique dans $X$ de sorte que 
$\tilde X~:=B_C(X)$ est projective,
\item [(ii)] il existe une 
courbe torique $\omega$
de $\tilde X$ Mori extr\'emale 
v\'erifiant $E\cdot\omega < 0$ o\`u $E$ est le diviseur exceptionnel
de $\pi~: \tilde X \to X$.
\end{enumerate}
Alors, $$ (C, N_{C/X}) \simeq (\PP ^1, 
{\mathcal O}_{\PP^1}(-1)^{\oplus n-1})$$
et il y a un diagramme 
\medskip

\centerline{
\xymatrix{  &\tilde{X}=B_C(X) \ar[dl]_{\pi} 
\ar[dr]^{\varphi _{\QQ ^+ [\omega]}}\\
X & & Y
}
}

\medskip 

o\`u $Y$ est une vari\'et\'e torique projective
contenant une sous-vari\'et\'e torique $Z$ 
v\'erifiant 
$$(Z, N_{Z/Y})\simeq (\PP ^{n-2}, 
{\mathcal O}_{\PP ^{n-2}}(-1)^{\oplus 2}),$$
$\varphi _{\QQ ^+ [\omega]}$
est la contraction extr\'emale d\'efinie par $\omega$ 
et l'\'eclatement de $Y$ le long de la sous-vari\'et\'e
$Z$ de sorte que $\varphi _{\QQ ^+ [\omega]}~: \tilde{X}\to Y$ 
et $\pi~: \tilde{X}\to  X$ ont m\^eme diviseur
exceptionnel $E$
v\'erifiant 
$$(E, N_{E/\tilde X})\simeq (\PP ^1\times
\PP ^{n-2}, {\mathcal O}_{\PP ^1 \times 
\PP ^{n-2}}(-1,-1)).$$

De plus, la classe d'une courbe rationnelle 
incluse dans $Z$ n'est pas ex\-tr\^e\-ma\-le dans le c\^one de
Mori des courbes effectives de $Y$.

\end{prop}

\noindent {\bf D\'emonstration.}
Comme $E \cdot \omega <0$, la courbe
$\omega$ est incluse dans $E$. Mais comme $X$ n'est pas projective,
la proposition~1 implique que $\omega$ n'est pas incluse dans une fibre
de $\pi$. Quitte \`a renum\'eroter les $e_i$, $1\leq i \leq n-1$,
on peut supposer que $\omega$ est d\'efinie par le mur
$<e_1,\ldots,e_{n-2},e>$ et on a une relation
$$ e_n + e_{n+1} + \sum_{i=1}^{n-2}a_i e_i + a e =0,$$
o\`u les $a_i$ et $a$ sont des entiers.
Le lemme suivant est l'analogue du lemme~2.

\begin{lemm} Les entiers $a_i$ sont n\'egatifs ou nuls
pour tout $i$ entre $1$ et $n-2$.
\end{lemm}  

\noindent {\bf D\'emonstration.}
Sinon, on a par exemple $a_1 >0$, et d'apr\`es Reid \cite{Rei83} 
Cor. 2.10(i),
les deux courbes 
$\omega ':= <e_2,\ldots,e_{n-2},e,e_n>$ et 
$\omega '':= <e_2,\ldots,e_{n-2},e,e_{n+1}>$
sont sur l'ar\^ete
extr\'emale $\QQ ^+ [\omega]$.
Or ces deux courbes sont contract\'ees par $\pi$, ce qui n'est pas
possible puisque $X$ est non projective.\finpreuve

\subsubsection*{Fin de la d\'emonstration de la proposition~5.}
Comme 
$$ -K_{\tilde X}\cdot \omega = 2 + \sum _{i=1}^{n-2}a_i+a >0$$
avec $a=E \cdot \omega <0$, il s'en suit que $a=-1$
et $a_i=0$ pour tout $1\leq i \leq n-2$.
La contraction extr\'emale $\varphi _{\QQ ^+ [\omega]}$
est donc divisorielle, de fibres non triviales 
isomorphes \`a $\PP ^1$, et comme
$$e_n+e_{n+1} = e = \sum _{i=1}^{n-1}e_i,$$
la courbe $C$ a son fibr\'e normal isomorphe
\`a ${\mathcal O}_{\PP^1}(-1)^{\oplus n-1}$, et donc 
$$(E, N_{E/\tilde X})\simeq (\PP ^1\times
\PP ^{n-2}, {\mathcal O}_{\PP ^1 \times 
\PP ^{n-2}}(-1,-1)).$$

%cette figure est exportee en 50,20
\begin{figure}[hbtp]
  \begin{center}
    \leavevmode
    \input{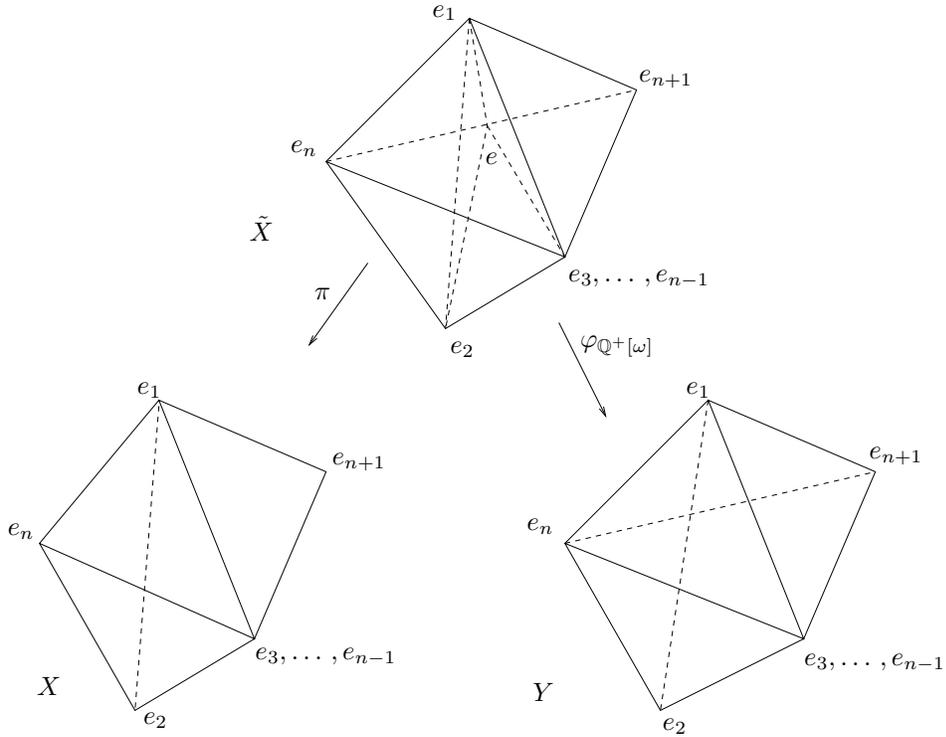}
    \caption{Flip interdit}
    %\label{fig:}
  \end{center}
\end{figure}

Il en r\'esulte - voir la figure~7 - 
que la contraction $\varphi_{\QQ ^+[\omega]}$
envoie $\tilde X$ sur la vari\'et\'e projective 
lisse $Y$
dont l'\'eventail est obtenu en rempla\c cant dans
$\tilde{\Delta}$ l'\'etoile de $e$ par la r\'eunion
des $n-1$ c\^ones maximaux 
$$<e_1,\ldots,\hat{e_i},\ldots,e_{n-1},e_n,e_{n+1}>,\,
1\leq i \leq n-1$$ et leurs faces.

Autrement dit, 
$\varphi_{\QQ ^+[\omega]}~: \tilde X \to Y$
est l'\'eclatement de $Y$ le long de la sous-vari\'et\'e torique
de codimension deux, isomorphe
\`a $\PP ^{n-2}$ et d\'efinie par le c\^one $<e_n,e_{n+1}>$.

Enfin, la classe d'une courbe rationnelle 
incluse dans $Z$ n'est pas ex\-tr\^e\-ma\-le dans le c\^one de
Mori des courbes effectives de $Y$ car $X$
n'est pas projective.
\finpreuve

\subsubsection*{D\'emonstration de la proposition~3.}
La proposition~3, {\S}3.2 pr\'ecise la proposition~5.
Si $X$ n'est pas projective, nous venons
de voir que la classe d'une courbe rationnelle 
incluse dans $Z$ n'est pas 
ex\-tr\^e\-ma\-le dans le c\^one de
Mori des courbes effectives de $\bar X$. R\'eciproquement, 
soit $\bar F$ une courbe rationnelle incluse dans $Z$, non 
ex\-tr\^e\-ma\-le dans le c\^one de
Mori des courbes effectives de $\bar X$~:
$[\bar F] = [C_1] + [C_2]$, o\`u les $[C_i]$
sont effectives, non incluses dans $Z$.
Si $C'_i$ d\'esigne la transform\'ee stricte de
$C_i$ pour $\bar {\pi}$, et si $F$ est une courbe de
$\tilde X$, incluse dans une fibre de $\pi$ telle
que $\bar{\pi}_* F = \bar F$, on a 
$$ [F] =  [C'_1] + [C'_2] + a [C'],$$
o\`u $a$ est un entier et $[C']$ la classe d'une fibre non
triviale de $\bar{\pi}$.
Or $$ -1 = E\cdot F= E\cdot C'_1+ E\cdot C'_2 + a E\cdot C'
= E\cdot C'_1+ E\cdot C'_2 - a,$$
d'o\`u 
$$ a = 1 + E\cdot C'_1+ E\cdot C'_2 \geq 1$$
car $E \cdot C'_i \geq 0$.
Ainsi, la classe de $[F]$ n'est pas extr\'emale dans 
le c\^one de
Mori des courbes effectives de $\tilde X$, donc
$X$ n'est pas projective d'apr\`es la proposition~1.
\finpreuve

\section{Quelques exemples}
Dans cette derni\`ere partie, nous donnons
des exemples explicites de vari\'et\'es toriques non projectives
pour toutes les situations pr\'ec\'edemment \'etudi\'ees.

\subsection{Exemples en dimension trois.}
\subsubsection{Flip interdit.}
L'exemple le plus c\'el\`ebre est donn\'e par la seule vari\'et\'e
torique $X$ non projective de dimension $3$
dont le groupe de Picard est de rang $4$~: son \'eventail est donn\'e
par la figure suivante projet\'ee
sur la face $<n_1,n_2,n_3>$ (voir \cite{Oda88} page~85),
o\`u $n_1$, $n_2$ et $n_3$ forment une base du r\'eseau $N$,
o\`u $n_0 = -n_1-n_2-n_3$ et pour $1\leq i \leq 3$, $n'_i=n_0+n_i$.

%cette figure est exportee en 50,20
\begin{figure}[hbtp]
  \begin{center}
    \leavevmode
    \input{fig11.pstex_t}
    \caption{}
    %\label{fig:}
  \end{center}
\end{figure}

Si $C_1$ (respectivement $C_2$, $C_3$) est la courbe de c\^one
$<n'_2,n'_3>$ (respectivement $<n'_1,n'_3>$, $<n'_1,n'_2>$),
alors $N_{C_i/X} = {\mathcal O}_{\PP ^1}(-1)^{\oplus 2}$,
la vari\'et\'e $\tilde {X_i} = B_{C_i}(X)$ est projective
et $X$ est obtenue {\em via}
un flip interdit d'une vari\'et\'e torique projective $Y_i$.

\subsubsection{Transformation \'el\'ementaire.}
D'apr\`es la classification \cite{Oda78} \cite{Oda88}
des vari\'et\'es toriques de dimension $3$
avec petit groupe
de Picard, les vari\'et\'es toriques non projectives de dimension
$3$ obtenues {\em via} une transformation \'el\'ementaire
d'une vari\'et\'e torique projective ont 
un groupe de Picard de rang sup\'erieur ou \'egal \`a $5$. L'\'eventail
le plus simple est celui de la vari\'et\'e $X_{a,b}$
suivante ($a$ et $b$ dans $\ZZ$) avec $\rho(X_{a,b})=5$ 
(labell\'ee $(4^4 5^4)'$
dans \cite{Oda88} page~64 et $[8-13']$
dans \cite{Oda78} page~79).

%cette figure est exportee en 50,20
\begin{figure}[hbtp]
  \begin{center}
    \leavevmode
    \input{fig12.pstex_t}
    \caption{}
    %\label{fig:}
  \end{center}
\end{figure}

\newpage 

La proposition suivante est imm\'ediate~:

\begin{prop} Pour tous $a$ et $b$ entiers,
les vari\'et\'es $X_{a-1,b}$, $X_{a+1,b}$, $X_{a,b-1}$
et $X_{a,b+1}$ sont obtenues {\em via} une transformation 
\'el\'ementaire de la vari\'et\'e $X_{a,b}$.
\end{prop}

\noindent {\bf D\'emonstration.} La figure suivante
est suffisante~:

%cette figure est exportee en 50,20
\begin{figure}[hbtp]
  \begin{center}
    \leavevmode
    \input{fig13.pstex_t}
    \caption{}
    %\label{fig:}
  \end{center}
\end{figure}

\finpreuve

\medskip

Le r\'esultat suivant est \'enonc\'e sans d\'emonstration
dans \cite{Oda78} page 80.

\begin{prop}
La vari\'et\'e $X_{a,b}$ est non projective
si et seulement si $a\neq 0$ et $b\neq -1$.
\end{prop}

\noindent {\bf D\'emonstration.} Si par exemple $a=0$
(le cas $b=-1$ se traite de fa\c con identique),
$X_{0,b}$ est projective
car obtenue
en \'eclatant un point puis une courbe \`a partir
d'une vari\'et\'e $Z$ projective (le groupe de Picard de $Z$
est de rang $3$ ou plus simplement $Z$
admet une fibration sur $\PP ^1$ de fibre 
$\PP ({\mathcal O}_{\PP^1} \oplus {\mathcal O}_{\PP^1}(b))$).

Supposons ensuite $a \geq 1$ et $b\geq 0$ et montrons
que $X_{a,b}$ est non projective (les trois autres cas se
traitent de fa\c con identique). Par l'absurde, si
$X_{a,b}$ est projective, il existe une fonction
$\varphi~: N_{\QQ} \to \RR$ lin\'eaire sur chaque c\^one
de l'\'eventail $\Delta$ d\'efinissant $X_{a,b}$ et strictement
convexe au sens o\`u $\varphi (n_1 +n_2) > \varphi (n_1) + \varphi (n_2)$
d\`es que $n_1$ et $n_2$ n'appartiennent pas \`a un m\^eme
c\^one de $\Delta$.
Ecrivons alors 
$$ (-n + bn') + (b+1) (-n') = -n -n' = n'' + (-n-n'-n'').$$
Comme $n''$ et $ -n-n'-n''$ n'appartiennent pas \`a un m\^eme c\^one
alors que $-n+bn'$ et $(b+1)(-n')$
appartiennent \`a $<-n+bn',-n'>$, on a~: 
$$ \varphi(-n + bn') +(b+1)\varphi(-n')> \varphi( n'') +\varphi(-n-n'-n'').$$
De m\^eme, $ (-n-n'-n'')+ (b+1)n'= -n''+(-n + bn')$ donne
$$ \varphi(-n-n'-n'')  +(b+1)\varphi(n')>\varphi( -n'')+ \varphi(-n + bn').$$
Par somme~:
$$(b+1)( \varphi(n')+\varphi(-n'))> \varphi( n'') +\varphi( -n'').$$
A nouveau, 
$$ (n+n'+an'') + a(-n'') = n+ n' \,\,  \mbox{ et } 
\,\, n+an''= (n+n'+an'')+(-n') $$
donnent
$$\varphi(n+n'+an'')+a\varphi(-n'') > \varphi( n) +\varphi( n') 
\,\,  \mbox{ et } \,\,
 \varphi( n) +a \varphi(n'')> \varphi(n+n'+an'')+\varphi( -n') $$
d'o\`u 
$$ a(\varphi( n'') +\varphi( -n''))> \varphi(n')+\varphi(-n').$$
En regroupant, 
$$ a(b+1)(\varphi( n'') +\varphi( -n''))>\varphi( n'') +\varphi( -n''),$$
et comme $0> \varphi( n'') +\varphi( -n'')$, on en d\'eduit $a(b+1) <0$,
ce qui est absurde.
\finpreuve

\medskip

Des deux propositions pr\'ec\'edentes se d\'eduit 
imm\'ediatement le~:

\medskip

\noindent {\bf Corollaire~D. } {\em Pour
$b\neq -1$, les vari\'et\'es toriques non projectives
$X_{1,b}$
et $X_{-1,b}$ sont obtenues {\em via} une transformation 
\'el\'ementaire de la vari\'et\'e torique projective $X_{0,b}$.
Pour $a\neq 0$, les vari\'et\'es toriques non projectives
$X_{a,0}$
et $X_{a,-2}$ sont obtenues {\em via} une transformation 
\'el\'ementaire de la vari\'et\'e torique projective $X_{a,-1}$. 
} 

\subsection{Des exemples en dimensions sup\'erieures~: 
une construction d'Ewald revisit\'ee.}

Dans cette derni\`ere partie, nous \'etudions g\'eom\'etriquement
et g\'en\'eralisons 
une cons\-truction
due \`a Ewald \cite{Ewa86} pour cons\-truire en toutes dimensions
sup\'erieures \`a quatre des vari\'et\'es toriques non projectives
avec petit groupe de Picard.

\subsubsection{Description g\'eom\'etrique, 
cas g\'en\'eral.}
Soit $X$ une vari\'et\'e analytique complexe 
(non n\'ecessairement torique) de dimension $n$
munie d'une action de $\CC ^*$.
On construit une vari\'et\'e $\tilde X$ de dimension $n+1$
en recollant $\CC \times X $ \`a $\PP ^1 \setminus \{0\} \times X$
sur l'ouvert commun $\CC ^* \times X$ gr\^ace
\`a l'action de $\CC ^*$~: $(t,x) \in \CC ^* \times X \subset \CC \times X $
est identifi\'e \`a 
$(t,t\cdot x) \in \CC ^* \times X \subset \PP ^1 \setminus \{0\} \times X$.
La vari\'et\'e $\tilde X$ poss\`ede par construction une fibration 
$f~: \tilde X \to \PP ^1$ dont 
toutes les fibres sont isomorphes \`a $X$.
Dans la suite, on notera $X_0 := f ^{-1}(0)$
et $X_{\infty} := f ^{-1}(\infty)$.
Remarquons qu'\`a toute sous-vari\'et\'e $Y$ de $X$, stable
par $\CC ^*$, 
correspond une sous-vari\'et\'e $\tilde Y$
de ${\tilde X}$ de m\^eme codimension que $Y$  
obtenue par le recollement de  
$\CC \times Y $ \`a $\PP ^1 \setminus \{0\} \times Y$.
Comme $\tilde X$, $\tilde Y$ poss\`ede une fibration sur $\PP ^1$,
de fibre isomorphe \`a $Y$. De plus,  
$$ \tilde Y \cap X_0 \simeq  \tilde Y \cap X_{\infty} \simeq Y.$$ 

\subsubsection{Description g\'eom\'etrique,
cas torique.}
Si $X$ est une vari\'et\'e torique, et si
$v$ est un \'el\'ement du r\'eseau $N$, on note classiquement $\lambda_v$
le sous-groupe \`a un param\`etre du tore $\Ho (M,\CC ^*)$
lui correspondant. Ce sous-groupe d\'efinit une 
action de $\CC ^*$ sur $X$ si bien que la construction
pr\'ec\'edente s'applique~: on note ${\tilde X}_v$
la vari\'et\'e correspondante et $f_v~: {\tilde X}_v \to \PP^1$
la fibration correspondante. Comme les actions sur $X$ de 
$\Ho (M,\CC ^*)$ et de $\lambda_v$ commutent, elles induisent
une action de $\Ho (M,\CC ^*)\times \CC^*$ sur ${\tilde X}_v$, faisant
de cette vari\'et\'e une vari\'et\'e torique de dimension $n+1$ 
poss\'edant une fibration $f_v$ \'equivariante lorsque
$\PP ^1$ est muni de l'action standard de $\CC ^*$.
Dans ce cas, $X_0$ et $X_{\infty}$
sont deux sous-vari\'et\'es toriques de ${\tilde X}_v$ et
si $Y$ est une sous-vari\'et\'e torique de $X$,
la sous-vari\'et\'e $\tilde Y$
de ${\tilde X}_v$ lui correspondant est
torique. 
De plus, il
y a dans ${\tilde X}_v$ une sous-vari\'et\'e torique privil\'egi\'ee,
que
l'on notera
${\tilde Y}_v$~: en effet $v$ appartient \`a l'int\'erieur
relatif d'un unique c\^one $\sigma _v$ de
l'\'eventail $\Delta$ et la sous-vari\'et\'e torique $Y_v$
de $X$ d\'etermin\'ee par $\sigma _v$ est fix\'ee par le sous-groupe
$\lambda _v$. Par cons\'equent, la sous-vari\'et\'e torique de 
${\tilde X}_v$ lui correspondant est isomorphe    
\`a $\PP ^1 \times Y_v$. 

\subsubsection{Description ``\'eventail''.}    
Si $\Delta$ est l'\'eventail dans $N_{\QQ}$ de $X$, 
l'\'eventail $\tilde{\Delta}_v $ 
dans $N_{\QQ} \oplus \QQ$ de $\tilde{X}_v$
est d\'efini de la fa\c con suivante~:
\begin{enumerate}
\item[-] on identifie $\Delta$ avec son image 
{\em via} l'inclusion $N_{\QQ}\oplus 0 \hookrightarrow N_{\QQ}\oplus \QQ$,
\item[-] les c\^ones de $\tilde{\Delta}_v $ sont exactement les 
$$ \tilde{\sigma}:=(\sigma,0) \,\,  , \,\,    
\tilde{\sigma}_+ :=(\sigma,0) + \QQ ^+ (v,1) \, \, \mbox{ et } \,\,
\tilde{\sigma}_- :=(\sigma,0) + \QQ ^+(0,-1)$$
o\`u $\sigma$ d\'ecrit l'ensemble des c\^ones de $\Delta$.
\end{enumerate}
 
\subsubsection{La construction d'Ewald.} 
Soit $X$ une vari\'et\'e torique de dimension $n$, 
$D$ un diviseur irr\'eductible torique de $X$. On note $v_D$
le g\'en\'erateur minimal dans $N$ du c\^one de dimension
un d\'efinissant $D$. D'apr\`es ce qui pr\'ec\`ede, 
la vari\'et\'e torique ${\tilde X}_{v_D}$ de dimension $n+1$
construite en $5.2.1-2$ contient un 
diviseur torique not\'e ${\tilde D}$ isomorphe \`a $\PP ^1 \times D$
et Ewald remarque que le fibr\'e normal ${\mathcal N}$
de ${\tilde D}$ dans ${\tilde X}_{v_D}$
satisfait 
$$ {\mathcal N}_{| \PP ^1 \times \{*\}} \simeq 
{\mathcal O}_{\PP ^1}(-1).$$
Par cons\'equent, il y a une vari\'et\'e torique 
$X_{ v_D}$ de dimension $n+1$ contenant
une sous-vari\'et\'e torique $Z_{ v_D}$ de codimension deux isomorphe 
\`a $D$ de sorte que ${\tilde X}_{v_D}$ soit l'\'eclatement de 
$X_{ v_D}$ le long de $Z_{ v_D}$, le diviseur exceptionnel 
de l'\'eclatement \'etant \'egal \`a ${\tilde D}$.
Remarquons que $X_{ v_D}$ contient deux sous-vari\'et\'es
toriques de codimension un isomorphes \`a
$X_0$ et $X_{\infty}$, dont l'intersection est isomorphe \`a 
$D$.

Le r\'esultat suivant est d\^u \`a Ewald~:

\medskip

\noindent {\bf Th\'eor\`eme.} {\em 
Soit $X$ une vari\'et\'e torique, $D$ un diviseur
torique de $X$
et $X_{ v_D}$ la vari\'et\'e torique construite pr\'ec\'edemment.
Alors $X_{ v_D}$ est projective si et seulement si
$X$ est projective. De plus, $X$ et $X_{ v_D}$ ont m\^eme
nombre de Picard.
}

\medskip

D'apr\`es ce qui pr\'ec\`ede, seule l'affirmation ``$X$ projective
implique $X_{ v_D}$ projective'' est non triviale, Ewald
la d\'emontre en consid\'erant l'\'eventail d\'efinissant $X_{ v_D}$. Ceci
peut aussi se d\'emontrer facilement en utilisant notre 
crit\`ere de projectivit\'e {\S}1.4. 

\subsubsection{Applications.}

Soit $X$ une vari\'et\'e torique non projective
de dimension $n$
devenant projective apr\`es \'eclatement le long d'une
sous-vari\'et\'e torique $Y$. On note $X' := B_Y(X)$.
Deux constructions diff\'erentes
sont possibles~:

\begin{enumerate}
\item[(i)] si $D$ est un diviseur irr\'eductible torique de
$X$ disjoint de $Y$, alors 
la vari\'et\'e torique $X_{ v_D}$ de dimension $n+1$
est non projective et devient projective 
apr\`es \'eclatement le long de
la sous-vari\'et\'e torique $\tilde Y$ correspondant \`a $Y$.
Par construction, la codimension de $Y$ dans $X$ est \'egale
\`a celle de $\tilde Y $ dans $X_{ v_D}$.
\item[(ii)] si $D$ est un diviseur irr\'eductible torique de
$X$ contenant $Y$, la transform\'ee stricte de $D$
dans $X'= B_Y(X)$ \'etant not\'ee $D'$, alors 
la vari\'et\'e torique $X_{ v_D}$ de dimension $n+1$
est non projective et devient projective apr\`es \'eclatement le long 
d'une sous-vari\'et\'e torique isomorphe \`a $Y$ (et que l'on 
note encore $Y$)
incluse dans $Z_{ v_D}$. En effet, on 
a $B_Y (X_{ v_D}) = X'_{v_{D'}}$ et cette derni\`ere vari\'et\'e est 
projective puisque $X'$ l'est.
Par construction, $X$ et $X_{ v_D}$ deviennent projectives apr\`es
\'eclatement le long de sous-vari\'et\'es isomorphes, donc de m\^eme
dimension !
\end{enumerate}

On en d\'eduit le corollaire suivant~:

\medskip

\noindent {\bf Corollaire~E. } {\em 
Pour tout entier $n \geq 3$, il existe une vari\'et\'e 
torique non projective, dont le nombre de Picard est 
\'egal \`a $4$, devenant projective apr\`es
\'eclatement le long d'une courbe torique.
}

\medskip

\noindent {\bf D\'emonstration.} 
Il suffit d'it\'erer la construction expliqu\'ee en (ii) 
\`a partir de la vari\'et\'e d\'ecrite en 5.1.1. \finpreuve

\medskip

Les vari\'et\'es du corollaire~E sont 
obtenues {\em via}
un flip interdit d'une vari\'et\'e torique projective, 
la m\^eme construction partant par exemple des vari\'et\'es $X_{1,b}$
($b\neq -1$) de 5.1.2 permet de construire, pour tout entier $n \geq 3$,
une vari\'et\'e 
torique non projective, dont le nombre de Picard est 
\'egal \`a $5$, obtenue
{\em via} une transformation \'el\'ementaire
d'une vari\'et\'e torique projective.

-----------

{\em 
\noindent Institut Fourier, UFR de Math\'ematiques\\
\noindent Universit\'e de Grenoble 1 \\
\noindent UMR 5582 \\
\noindent BP 74 \\
\noindent 38402 SAINT-MARTIN d'H\`ERES \\
\noindent FRANCE \\
\noindent e-mail : bonavero@ujf-grenoble.fr
}
\end{document}